\documentclass{aims}
\usepackage{amsmath}
  \usepackage{paralist}
  \usepackage{graphics} 
  \usepackage{epsfig} 
\usepackage{graphicx}  \usepackage{epstopdf}
 \usepackage[colorlinks=true]{hyperref}
\hypersetup{urlcolor=blue, citecolor=red}

\def\currentvolume{X}
 \def\currentissue{X}
  \def\currentyear{200X}
   \def\currentmonth{XX}
    \def\ppages{X--XX}
     \def\DOI{10.3934/xx.xx.xx.xx}

\newtheorem{theorem}{Theorem}[section]

\newtheorem{lemma}[theorem]{Lemma}

\theoremstyle{definition}
\newtheorem{definition}[theorem]{Definition}
\newtheorem{remark}{Remark}

\newcommand{\ep}{\varepsilon}
\newcommand{\eps}[1]{{#1}_{\varepsilon}}
\newcommand{\la}{\lambda}

\newcommand{\mc}[1]{\mathcal{#1}}
\newcommand{\comment}[1]{}

\newcommand{\mdm}[1]{\mathrm{#1}}
\newcommand{\md}{\mathrm{d}}

\newfont{\smoldita}{cmmib8}
\newfont{\boldita}{cmmib10}
\newfont{\bboldita}{cmmib10 scaled\magstep1}
\newcommand{\sem}[1]{\mbox{$({#1}(t))_{t \geq 0}$}}

\newcommand{\mbb}[1]{\mathbb{#1}}

\newcommand{\nn}{\nonumber}
\newcommand{\e}{\epsilon}

\newcommand {\ti}[1]{\widetilde {#1}}

\newcommand{\p}{\partial}
\newcommand{\cl}[2]{\int\limits_{#1}^{#2}}

\newcommand{\fin}[1]{{#1}^{in}}

\title[Global solutions of continuous C-F equations] 
      {Global solutions of continuous coagulation--fragmentation equations with unbounded coefficients}

\author[Jacek Banasiak]{}


 \email{jacek.banasiak@up.ac.za}

\thanks{The research has been partially supported by the National
Science Centre of Poland Grant 2017/25/B/ST1/00051 and the National Research Foundation of South Africa Grant 82770}

\begin{document}
\maketitle

\centerline{\scshape Jacek Banasiak}
\medskip
{\footnotesize
 \centerline{Department of Mathematics and Applied Mathematics}
   \centerline{University of Pretoria}
   \centerline{ Pretoria, South Africa}
} 

\bigskip
\centerline{\small The paper is dedicated to Gis\'ele Ruiz Goldstein on the occasion of her birthday}
\bigskip
 \centerline{(Communicated by the associate editor name)}

\begin{abstract}
In this paper we prove the existence of global classical solutions to continuous coagulation--fragmentation equations with unbounded coefficients under the sole assumption that the coagulation rate is dominated by a power of the fragmentation rate, thus improving upon a number of recent results by not requiring any polynomial growth bound for either rate. This is achieved by proving a new result on the analyticity of the fragmentation semigroup and then using its regularizing properties to prove the local and then, under a stronger assumption, the global classical solvability of the coagulation--fragmentation equation considered as a semilinear perturbation of the linear fragmentation equation.   Furthermore, we show that weak solutions of the coagulation--fragmentation equation, obtained by the weak compactness method, coincide with the classical local in time solutions provided the latter exist.
\end{abstract}

\section{Introduction}

Coagulation equations, introduced by Smoluchowski \cite{Smo16,Smo17} in the discrete case and in \cite{muller1928allgemeinen} in the continuous one,  and extended in \cite{becker, blatz1945, Melz57b, McZi87, vigi} to include the reverse fragmentation processes, have proved crucial in numerous applications, ranging from polymerization, aerosol formation, animal groupings, phytoplankton dynamics, to rock crushing and planetesimals formation, see a survey in \cite[Vol. I]{BLL} and, as such, they have been extensively studied in engineering, physical and mathematical sciences. We note that coagulation--fragmentation processes can be studied also in a probabilistic setting, see e.g. \cite{Bert06}, but in this paper we shall focus on the deterministic approach that yields the following kinetic type
 continuous coagulation--fragmentation equation
\begin{subequations}\label{PhPr001}
\begin{equation}
\partial_t f(t,x) = \mathcal{C}f(t,x) + \mathcal{F}f(t,x)\ , \qquad  (t,x)\in \mbb R_+^2, \label{PhPr001a}
\end{equation}
with the initial condition
\begin{equation}
f(0,x) = f^{in}(x)\ , \qquad x\in \mbb R_+,  \label{PhPr001b}
\end{equation}
\end{subequations}
where  the coagulation operator $\mathcal{C}$ and the fragmentation operator  $\mathcal{F}$ operators are,  respectively, given by
\begin{equation}
\mathcal{C}f(x) = \frac{1}{2} \int_0^x k(x-y,y) f(x-y) f(y)\ \mathrm{d}y - \int_0^\infty k(x,y) f(y) f(x)\ \mathrm{d}y \label{PhPr002}
\end{equation}
and
\begin{equation}
\mathcal{F}f(x) = \mc A f(x)+ \mc B f(x) = - a(x) f(x) + \int_x^\infty a(y) b(x,y) f(y)\ \mathrm{d}y, \label{PhPr003}
\end{equation}
for $x\in \mbb R_+$. In \eqref{PhPr001}, \eqref{PhPr002} and \eqref{PhPr003}, $f$ is the density of particles of mass $x$, the coagulation kernel $k$ is a nonnegative and measurable symmetric function defined on $\mbb R_+^2$, the overall fragmentation rate $a$ is a nonnegative measurable function on $\mbb R_+$ satisfying \begin{equation}
a \in L_{\infty,loc}([0,\infty)).
\label{aloc}
\end{equation}
The daughter distribution function $b,$ sometimes referred to as the fragmentation kernel, is a nonnegative and measurable function such that for a.e. $y>0$,
\begin{equation}
\int_0^y x b(x,y)\ \mathrm{d}x = y \;\;\text{ and }\;\; b(x,y) = 0 \;\text{ for a.e. }\; x>y. \label{PhPr004}
\end{equation}
We recall that the first condition in \eqref{PhPr004} ensures that there is no loss of matter during fragmentation events.

Typically, the analysis of the coagulation--fragmentation equations  is done in the spaces $X_0 := L_1(\mathbb{R} _+, \md x)$ and $X_1:=L_1(\mathbb{R} _+,  \md x)$ since the norm of a nonnegative $f$ in $X_0$,
$$
\|f\|_{[0]} = \cl{0}{\infty} f(x)\md x,
$$
gives the total number of particles in the system, while
$$
\|f\|_{[1]} = \cl{0}{\infty} f(x)x\md x
$$
gives its total mass. It follows that by introducing some control on the evolution of large particles we can improve the properties of the involved equations.  The easiest way to introduce such a control is to consider the problem in the spaces $X_m := L_1(\mathbb{R} _+, x^m \md x)$ { and } $X_{0,m} := L_1(\mathbb{R}_+,(1+x^m) \md x)$; the natural norms in these spaces will be denoted by $\|\cdot\|_{[m]}$ and $\|\cdot\|_{[0,m]}$.  To shorten notation, we use the abbreviation $w_m(x) =1+x^m$.

 In its full generality, \eqref{PhPr001} is a nonlinear, nonlocal integro-differential equation with unbounded coefficients and hence its solvability presents a wide array of challenges. Early attempts, whose survey can be found in \cite{BLL}, mostly focused on finding particular solutions to \eqref{PhPr001}, often by quite ingenious methods. Systematic mathematical studies  of \eqref{PhPr001} date back to the 1980s and two main (deterministic) ways to approach have emerged. In the so called weak compactness method, used originally in e.g. \cite{BaCa90} for the discrete version of the problem and in  \cite{Stew89} for the continuous one (see also \cite{ELMP03} for a more comprehensive approach), the equation is first truncated to  yield a  more tractable family of problems approximating \eqref{PhPr001}. Then it is shown that the solutions to the truncated problems form a weakly compact  family of functions from which one can select a subsequence converging in a suitable topology to a solution of an appropriate weak formulation of \eqref{PhPr001}. The other method, which can be termed the operator one, was initiated in \cite{AizBak}, and, roughly speaking, consists in considering the coagulation part as a perturbation of the linear fragmentation part. Then the theory of semigroups of operators is used to first obtain the (linear) fragmentation semigroup and hence  solve \eqref{PhPr001} by an appropriate fixed point technique. Each method has its advantages and disadvantages that make them suitable for different scenarios and thus they have been developed to large extent independently of each other. The weak compactness approach mostly uses the properties of the coagulation part and can deliver the existence of a solution for a large class of coagulation kernels but then the fragmentation part must somehow match the coagulation term; also other properties, such as mass conservation, regularity, or uniqueness of solutions, have to be proved independently under much more stringent assumptions. On the other hand, the operator method provides the existence of unique, mass conserving and classical solutions but, while being able to deal with even very singular fragmentation processes,  its applications to the full problem \eqref{PhPr001} for a long time were restricted to bounded coagulation kernels. This has changed in the recent few years with the realization that the fragmentation semigroup is analytic for a large class of fragmentation rates $a$ and the daughter distribution functions $b$. This, in turn, allowed for proving the  classical solvability of \eqref{PhPr001} even if the coagulation rate $k$ is unbounded as long as it is dominated in a suitable sense by the fragmentation rate $a$, see \cite{BaLa12a, BLL13} and \cite[Section 8.1.2]{BLL}. The proofs use interpolation spaces between $X_{0,m}$ and the domain of the generator of the fragmentation semigroup in this space.

The main aim of this paper is twofold. First, we extend the results of \cite{ BLL, BaLa12a} by removing the assumption that the fragmentation rate is of polynomial growth. This requires a new proof of the analyticity of the fragmentation semigroup that this time is based on the Miyadera-Voigt perturbation theorem, see e.g. \cite{Voi77} or \cite[Corollary III.3.16]{EN}, with the help of \cite[Lemma 4.15]{BaAr}.  It turns out that the local in time  solvability of \eqref{PhPr001} remains the same as in \cite{BaLa12a}, but the global one requires some new moment estimates: for the zeroth moment we adapt the ideas present in \cite{ELMP03, Stew91} (see also
\cite[Lemma 8.2.27]{BLL}), while for the estimates in the interpolation spaces we use the Henry-Gronwall inequality as in \cite{Banasiak2019}. We emphasize that, in contrast to e.g. \cite{DuSt96b, EMP02, ELMP03}, we do not require any polynomial estimates on the coagulation kernel, or on the fragmentation rate; also we allow the expected number of daughter particles to be unbounded for large size of the parent particle.   Second, we show that if the coefficients of \eqref{PhPr001} satisfy the assumptions of the local existence theorem, then the solutions to the truncated problems, constructed in e.g. \cite{ELMP03} as the approximations to a weak solution to \eqref{PhPr001} in the weak compactness method,  converge strongly to the classical solution of \eqref{PhPr001} on its maximal interval of existence, confirming thus the  fact, not entirely surprising,  that both methods agree with each other whenever they are both applicable.

\textbf{Acknowledgement.} The author is grateful to Prof. Mustapha Mokhtar-Kharroubi for the suggestions concerning the application of the Miyadera theorem to the problem and to Dr. S. Shindin, whose ideas for the estimates in the discrete case have helped to develop their counterparts in the continuous case.

\section{Analytic fragmentation semigroup}

Let $X$, $Y$ be Banach spaces. The space of bounded linear operators from $X$ into $Y$ will be denoted by $\mc L(X, Y),$ shortened to $\mc L(X)$ if $X=Y$.  If an operator $O$ generates a $C_0$-semigroup, this semigroup will be denoted by \sem{G_O}.

 For $m \geq 0$ we introduce
\begin{align}
n_m(y)&=\cl{0}{y}b(x,y)x^m\md x, \label{nmy}\\
N_m(y)& = y^m-n_m(y);
\label{Nmy}
\end{align}
then we have the inequalities,  \cite[Eqn. (2.3.16)]{BLL},
\begin{equation}
N_m(y) > 0, \quad m>1, \qquad N_1(y) = 0, \qquad N_m(y) <0, \quad 0\leq m<1.
\label{Nm}
\end{equation}
First, let us assume $m\geq 1$. We define $A_m f := -af$ on
$$
D(A_m) =\{f \in  X_{m}\; : \;af\in  X_{m}\}
$$
and, using the definition  (\ref{PhPr003}) of $\mc B$, by (\ref{nmy}) we easily get
\begin{equation}
\|\mc{B}f\|_{[m]} = \cl{0}{\infty}a(y)n_m(y)f(y) \mdm{d}y <\infty, \quad f\in D(A_m)_+. \label{expB}
\end{equation}
 Hence, we can define $B_m = \mc B|_{D(A_m)}$.

Moving to $X_{0,m}$ we need to introduce some control on the number of particles produced in fragmentation events.  Hence, besides (\ref{aloc}) and (\ref{PhPr004}),  we assume that there is $l\geq 0$ and $b_0\in \mbb R_+$ such that for any $x\in \mbb R_+$
\begin{equation}
n_0(x)=\cl{0}{x}b(y,x) \mdm{d}y \leq b_0(1+x^l).
 \label{ass2}
\end{equation}
Similarly to $A_m$, for any $m\geq 1$ we define $A_{0,m} f := -af$ on
$$
D(A_{0,m}) =\{f \in  X_{0,m}\; : \; af\in  X_{0,m}\}.
$$
Defining $B_{0,m}$ is, however, slightly more involved.
\begin{lemma}\label{lem34}
If $0\leq f\in D(A_{0,m})$ with $m\geq l$, then
\begin{equation}
\|\mc{B}f\|_{[0,m]} = \cl{0}{\infty}a(y)(n_m(y)+n_0(y))f(y) \mdm{d}y <\infty. \label{expB'}
\end{equation}
 \label{t3.4a}
\end{lemma}
\proof  Let $f\in D(A_{0,m})_+$. By (\ref{expB}), it suffices to estimate
\begin{align*}
&\cl{0}{\infty}\left(\cl{x}{\infty}a(y)b(x,y)f(y) \mdm{d}y\right)\!\! \mdm{d}x
=\cl{0}{\infty}a(y)f(y)\left(\cl{0}{y}b(x,y) \mdm{d}x\right) \mdm{d}y\label{diss''}\\&=\cl{0}{\infty}a(y)n_0(y)f(y) \mdm{d}y
\leq
b_0\cl{0}{\infty}a(y)(1+y^l)f(y) \mdm{d}y \leq 2 b_0\cl{0}{\infty}a(y)w_m(y) f(y) \mdm{d}y<\infty,\nn
\end{align*}
where we used the estimate
\begin{equation}
1+y^l\leq 2(1+ y^m), \label{4}
\end{equation}
 if $m\geq l$.
Hence, we can define $B_{0,m} = \mc B|_{D(A_{0,m})}$  provided $m\geq l$. \hfill\qed

\begin{theorem} \label{thnewchar0} Let $a,b$ satisfy  (\ref{aloc}),  (\ref{PhPr004}) and (\ref{ass2}).  Let further
for some $m_0>1$
\begin{equation}
\liminf\limits_{y\to \infty}\frac{N_{m_0}(y)}{y^{m_0}} >0.
\label{goodchar}
\end{equation}
Then
\begin{enumerate}
 \item (\ref{goodchar}) holds for all $m >1$;
 \item $F_{0,m} := A_{0,m}+B_{0,m}$ is the generator of a positive analytic semigroup  \sem{G_{F_{0,m}}}, on $X_{0,m}$ for any $m>\max\{1,l\}$.
     \end{enumerate}
\end{theorem}
\proof 1. This result in the discrete case is due to \cite{Banasiak2019}. Let $y\geq 1$ and $m_0>1$. It is easy to see that \eqref{goodchar} is equivalent to the existence of a constant $\delta_{m_0}>0$ such that $\inf_{y\geq 1}N_{m_0}/y^{m_0} \geq \delta_{m_0}$.  We have
$$
\frac{d}{dm} \frac{N_m(y)}{y^m} =    - \frac{1}{y^m}\cl{0}{y} b(x,y)x^m\ln\left(\frac{x}{y}\right )\md x >0
$$
and
$$
\frac{d^2}{dm^2} \frac{N_m(y)}{y^m} =     -\frac{1}{y^m}\cl{0}{y} b(x,y)x^m\ln^2\left(\frac{x}{y}\right)\md x <0,
$$
 where the differentiation under the sign of the integral is justified as $x^{m-1}(\ln x)^i$, $i=1,2$, is bounded due to $m>1$. Hence, if
 $\frac{N_{m_0}(y)}{y^{m_0}} > \delta_{m_0}$ for some $\delta_{m_0}> 0$ and some $m_0>1$, then $\frac{N_m(y)}{y^m}>\delta_{m_0}$ for any $m> m_0$. Further, the inequality for the second derivative shows that $m \mapsto \frac{N_m(y)}{y^m}$ is concave; that is, since $\frac{N_1(y)}{y} =0$, for $m\in (1,m_0]$ and $y\geq 1$ we obtain
 $$
 \frac{N_m(y)}{y^m} \geq \frac{N_{m_0}(y)}{y^{m_0}(m_0-1)}(m-1)\geq \frac{\delta_{m_0}(m-1)}{m_0-1},
 $$
 which gives \eqref{goodchar} in the interval $(1, m_0]$.

 2. To prove the generation result, we use the Miyadera--Voigt theorem, see \cite[Corollary III.3.16]{EN} and \cite[Lemma 4.15]{BaAr}. For this we observe that  (\ref{goodchar}) and the positivity of $m-l$ imply that there is $r>0$ such that for $x\geq r$ we have
\begin{equation}
\frac{n_m(x)}{x^m}  \leq c'<1, \qquad \frac{b_0(1+x^l)}{1+x^m} \leq  \frac{1-c'}{4},\label{pp}
\end{equation}
see (\ref{ass2}). Furthermore, by \eqref{ass2}, there is an $\zeta>0$ such that
\begin{equation}
\mathrm{ess}\!\!\!\sup\limits_{0\leq x\leq r} \frac{a(x)b_0(1+x^l)}{a(x)+\zeta}  \leq  \frac{1-c'}{4}.
\label{esss}
\end{equation}
Consider the operator $(A_{0,m}-\zeta I, D(A_{0,m})$. Then for $f \in D(A_{0,m})_+$ we have
\begin{align*}
&\int_0^\delta \|B_{0,m} G_{A_{0,m}-\zeta I}(t)f\|_{[0,m]}\md t \\
&= \int_0^\delta\left( \int_0^\infty (1+x^m) \left(\int_x^\infty a(y)b(x,y) e^{-(a(y)+\zeta)t} f(y)\md y\right)\md x\right)\md t\\
&\leq  \int_0^\infty (1+x^m) \left(\int_x^\infty \frac{a(y)b(x,y) }{a(y)+\zeta} f(y)\md y\right)\md x \\
&= \int_0^\infty \frac{a(y)f(y)}{a(y)+\zeta}\left(\int_0^y (1+x^m)b(x,y)\md x\right)\md y = I_1+I_2,
\end{align*}
where, by \eqref{nmy}, \eqref{ass2}, \eqref{esss} and the monotonicity of $x\mapsto 1+x^m,$
\begin{align*}
I_1 &:= \int_0^{r}\frac{a(y) f(y)}{a(y)+\zeta} \left(\int_0^y (1+x^m)b(x,y)\md x\right)\md y\leq \int_0^{r}\frac{a(y)(1+y^m)n_0(y)}{a(y)+\zeta} f(y) \md y\\
&\leq b_0\int_0^{r}\frac{a(y)(1+y^l)(1+y^m)}{a(y)+\zeta} f(y) \md y \leq \frac{1-c'}{4}\|f\|_{[0,m]}
\end{align*}
and, by \eqref{pp} and \eqref{4},
\begin{align*}
I_2&:= \int_{r}^\infty \frac{a(y)f(y)}{a(y)+\zeta}  \left(\int_0^y (1+x^m)b(x,y)\md x\right)\md y\\
&\leq \int_{r}^\infty \frac{a(y)n_0(y)f(y)}{a(y)+\zeta}  \md y + \int_{r}^\infty \frac{a(y)n_m(y) f(y)}{a(y)+\zeta} \md y\\
&\leq \|f\|_{[0,m]} \left(\mathrm {ess} \sup\limits_{y\geq r}  \frac{a(y)b_0(1+y^l)}{(1+y^m)(a(y)+\zeta)} + c' \right)\leq \frac{3c'+1}{4}\|f\|_{[0,m]}.
\end{align*}
Hence
$$
\int_0^\delta \|B_{0,m} G_{A_{0,m}-\zeta I}(t)f\|_{[0,m]}\md t \leq I_1+I_2 \leq \gamma \|f\|_{[0,m]}
$$
with $\gamma = (c'+1)/2 <1$. Therefore $B_{0,m}$ is a Miyadera perturbation of $A_{0,m}-\zeta I,$ and hence of $A_{0,m}$, see \cite[Lemma 4.15]{BaAr}. Using \cite[Exercise III.3.17(1)]{EN} and Arendt--Rhandi theorem, \cite{AR}, we conclude that $F_{0,m} = A_{0,m}+B_{0,m}$ is the generator of an analytic positive semigroup.
\section{Local solvability}
As mentioned in the introduction, the local in time solvability of \eqref{PhPr001} can be proved exactly as in \cite{BaLa12a}, see also \cite[Theorem 8.1.2.1]{BLL}. Certain notation and intermediate estimates will be, however, used in the proof of the global existence and thus are recalled below.

 We assume that $a$ and $b$ satisfy \eqref{aloc}, \eqref{PhPr004}, \eqref{ass2} and \eqref{goodchar}. Hence the fragmentation operator $(F_{0,m}, D( A_{0,m})) = (A_{0,m} + B_{0,m}, D(A_{0,m}))$ is the generator of an analytic positive semigroup on $X_{0,m}$ whenever $m>\max\{1,l\}.$ The coagulation kernel $k$ is assumed to be a measurable symmetric function such that,
 for some $K>0$ and $0 < \alpha<1$,
\begin{equation}
0\leq k(x,y) \leq K(1+a(x))^\alpha(1+a(y))^\alpha,  \quad (x,y) \in \mbb R_+^2.\label{kass1}
\end{equation}
 This assumption is sufficient for the local-in-time solvability of (\ref{PhPr001}). However, to prove that the solutions are global in time we need to strengthen (\ref{kass1}) to
\begin{equation}
0\leq k(x,y) \leq K\big((1+a(x))^\alpha+(1+a(y))^\alpha\big), \quad (x,y) \in \mbb R_+^2,
\label{kass2}
\end{equation}
again for $K > 0$ and  $0 < \alpha<1$.
Thus, using the linear operators $A_{0,m}$ and $B_{0,m}$, and the nonlinear operator $C_{0,m}$, defined via  (\ref{PhPr002}) but now only for $f$ in the maximal domain
\[
D(C_{0,m}) := \{f \in X_{0,m} : \mc C f \in X_{0,m}\},
\]
the initial-value problem  \eqref{PhPr001} can be written as the following abstract semilinear Cauchy problem in $X_{0,m}$:
\begin{equation}
\p_tf = A_{0,m}f + B_{0,m} f + C_{0,m}f, \qquad f(0) = f^{in}.
\label{feco1}
\end{equation}
We note that, in general, $0 \notin \rho(F_{0,m})$ and therefore to enable us to define appropriate intermediate spaces we consider
$$
F_{0,m,\omega}:= F_{0,m} - \omega I = A_{0,m} - \omega I + B_{0,m} = A_{0,m,\omega} + B_{0,m},
$$ where
$$
A_{0,m,\omega} : = A_{0,m} - \omega I,
$$
assuming that $\omega$ is greater than the type of \sem{G_{F_{0,m}}}.  We also assume that $\omega >1$ to simplify using (\ref{kass1}) and (\ref{kass2}). Clearly, $(F_{0,m,\omega}, D(A_{0,m}))$ is also the generator of an analytic semigroup $(G_{F_{0,m,\omega}}(t))_{t \ge 0} =  (e^{-\omega t}G_{F_{0,m}}(t))_{t \ge 0},$  but now we have the desired property that $0 \in \rho(F_{0,m,\omega})$.   Thus, as in \cite{BaLa12a},  the intermediate spaces between $D(F_{0,m,\omega}) = D(A_{0,m,\omega})$, see \cite[Section 2.2]{Lun}, satisfy
$$
D_{F_{0,m,\omega}}(\alpha, 1) = D_{A_{0,m,\omega}}(\alpha, 1) = X_{0,m}^{\alpha}, \qquad \alpha \in (0,1),
$$
where
\begin{equation}
 X_{0,m}^{\alpha} := \left\{ f \in X_{0,m}:\; \ \int_0^\infty |f(x)|(\omega + a(x))^\alpha w_m(x)\,\md x < \infty \right\},
 \label{frps1}
\end{equation}
and equality of the spaces is interpreted in terms of equivalent norms, see also the Stein--Weiss theorem \cite[Corollary 5.5.4]{bergh1976}.
 The natural norm on $X_{0,m}^{\alpha}$ will be denoted by $\|\cdot\|^{(\alpha)}_{[0,m]}$, and we note that  $X_{0,m}^0 =  X_{0,m}$, and $X_{0,m}^1 = D(A_{0,m,\omega})$.

 Hence, in general,  if \sem{G} is an analytic semigroup in $X_{0,m}$ satisfying $$\|G(t)\|_{\mc L(X_{0,m})} \leq M_{0,m}^{(0)}e^{\omega_{0,m}t},$$ then, by \cite[Theorem II.4.6(c)]{EN}, it is a family of operators  in $\mc L(X_{0,m}, X_{0,m}^{1})$ satisfying
 $$\|G(t)\|_{\mc L(X_{0,m},X^1_{0,m})} \leq M_{0,m}^{(1)}e^{\omega_{0,m}t}t^{-1},\quad t>0.$$ Then it follows that $G(t)\in \mc L(X_{0,m}, X_{0,m}^\alpha)$ and there is $M^{(\alpha)}_{0,m}$ such that
\begin{equation}\label{eq2.4a}
\|G(t)\|_{\mc L(X_{0,m},X^\alpha_{0,m})} \leq \frac{M_{0,m}^{(
\alpha)}e^{\omega_{0,m}t}}{t^\alpha}, \quad t>0,\, 0<\alpha<1.
\end{equation}

\begin{theorem} Assume that $a$ and $b$  satisfy \eqref{aloc}, \eqref{PhPr004}, \eqref{ass2}, \eqref{goodchar} and  let $m > \max\{1,l\}$.
Further, let $k$ satisfy (\ref{kass1}). Then, for each $f^{in}\in X_{0,m,+}^{\alpha}$, there is
$\tau(f^{in})>0$ such that the initial-value problem (\ref{feco1}) has a
unique nonnegative classical solution
\begin{equation}
f \in C\left([0,\tau(f^{in})),
X_{0,m}^{\alpha}\right)\cap C^1\left((0,\tau(f^{in})), X^{\alpha}_{0,m}\right)\cap C\left((0,\tau(f^{in})), D(A_{0,m})\right).
\label{fprop}
\end{equation}
 \label{fcth1}
\end{theorem}
\begin{remark} As we mentioned, the proof of this result is the same as of \cite[Theorem 2.2]{BaLa12a} which was proved under the additional assumption that $a$ is polynomially bounded. This assumption, however, was only needed to prove, by an alternative method, that \sem{G_{F_{0,m}}} is an analytic semigroup generated by $F_{0,m} = A_{0,m}+B_{0,m}$ in $X_{0,m}$ with a suitably bigger $m$ depending also on the growth rate of $a$. The only other consequence of the generation theorem of \cite{BaLa12a} is that \sem{G_{F_{0,m}}} is a quasi-contractive semigroup (that is, satisfying $\|G_{F_{0,m}}\|_{\mc L(X_{0,m})} \leq e^{\omega_m t}$ for some $\omega_m$) which in turn allowed in the proof of \cite[Theorem 2.2]{BaLa12a} to use the Trotter--Kato  representation formula to prove that certain auxiliary semigroups are positive. This result is, however, available also by a direct analysis of the construction of these semigroups.
\end{remark}
We recall some equalities and  inequalities used in the proof of \cite[Theorem 2.2]{BaLa12a} that will be used in the sequel. First, let $\ti{\mc C}$ denote the bilinear form obtained from $\mc C$; that is,
$$
\ti{\mc C}(f,f) = \mc Cf
$$
where $\mc C$ is defined in \eqref{PhPr002}. Then direct calculations yield,
for any measurable $\theta,$
\begin{eqnarray}
\int_0^\infty \,\theta(x)\,[\tilde{\mc{C}}(f,g)](x)\,\md x &=& \frac{1}{2}\int_0^\infty \int_0^\infty \theta(x+y)k(x,y)f(x)g(y)\,\md x\md y\nn\\
&&\phantom{xxx}- \int_0^\infty \int_0^\infty\theta(x)k(x,y) f(x)g(y)\,\md x\md y
\label{moments1}
\end{eqnarray}
and, by symmetry,
\begin{equation}
\int_0^\infty \,\theta(x)\,[\mc{C}f](x)\,\md x
= \frac{1}{2}\int_0^\infty\int_0^\infty \chi_{\theta}(x,y)k(x,y) f(x)f(y)\,\md x\,\md y\,,
\label{moments2}
\end{equation}
where
\[
\chi_{\theta}(x,y) = \theta(x+y) - \theta(x) - \theta(y).
\]
In particular, for $\theta(x) = 1 + x^m$ we will be using the elementary inequality
$$
(x+y)^m \leq 2^m(x^m + y^m), \qquad x,y\in \mbb R_+^2,\, m\geq 0,
$$
as well as
\begin{equation}
0 \le (x+y)^m - x^m - y^m \le c_m \left( x y^{m-1} + x^{m-1} y \right)\ , \qquad x,y\in \mbb R_+^2,\,m>1,\label{momestb}
\end{equation}
for some $c_m$, see \cite[Lemma 7.4.2]{BLL}.

Next,  \eqref{moments1} and (\ref{kass1}) with $f,g\in X_{0,m}^{\alpha}$ and $\theta = w_m(x) = 1+x^m$ yield
\begin{equation}
\|\tilde{\mc{C}}(f,g)\|_{[0,m]} \leq (1+2^m) K\|f\|_{[0,m]}^{(\alpha)}\|g\|_{[0,m]}^{(\alpha)}.\label{cest2}
\end{equation}
On the other hand, assumption \eqref{kass2} yields in a similar way
\begin{equation}\label{cest3}
\|\mc C f\|_{[0,m]} \leq 2^{m+1} K\|f\|_{[0,m]}^{(\alpha)} \|f\|_{[0,m]}.
\end{equation}

\section{Relation with weak solutions}

Weak solutions to \eqref{PhPr001} are constructed as weak limits of solutions $f_r$ to the problem \eqref{PhPr001} with the coefficients $a$ and $k$ modified as follows
\begin{equation}
a_r(x) = \left\{\begin{array}{lcl} a(x)&\text{for}& x\leq r\\0&\text{for}& x>r,\end{array}\right.\qquad k_r(x,y) = \left\{\begin{array}{lcl} k(x,y)&\text{for}& x+y\leq r\\0&\text{for}& x+y>r,
\end{array}\right.
\label{arkr}
\end{equation}
see e.g. \cite{ELMP03}, or \cite[Lemma 8.2.24]{BLL}.
\begin{theorem} Assume that the assumptions of Theorem \ref{fcth1} are satisfied and  $f$ is the solution to \eqref{feco1} satisfying \eqref{fprop}. If $(f_r)_{r>0}$ are approximate solutions defined above, then
\begin{equation}
\lim\limits_{r\to \infty} f_r = f
\label{approx}
\end{equation}
in $C([0,T], X^{\alpha}_{0,m})$ for any $T<\tau(f^{in})$.\label{thappr}
\end{theorem}
\proof Let $0<r<\infty$. We observe that the solutions $f_r, r>0,$ to the truncated problem \eqref{PhPr001} are unique and thus coincide with the solutions obtained by the semigroup method as follows.  Let us denote by $\mc F_r$ the fragmentation expression \eqref{PhPr003} restricted to $[0,r).$ Taking into account the fact  that the spaces $X_{0,m,r} = L_1((0,r), w_m(x)\md x)$ are invariant under the action of $\sem{G_{F_{0,m}}}$ we find that the (uniformly continuous and analytic) semigroups \sem{G_{F^{(r)}_{0,m}}} generated by $F^{(r)}_{0,m} = F_{0,m}|_{X_{0,m,r}}$ coincide with $\sem{G_{F_{0,m}}|_{X_{0,m,r}}}$   and thus all the estimates for \sem{G_{F^{(r)}_{0,m}}} are the same as for \sem{G_{F_{0,m}}} and independent of $r$. Also, \eqref{kass1} is satisfied for $a_r$ and $k_r$ uniformly in $r$ and thus the solution $f_r$ satisfies all estimates of Theorem \ref{fcth1}  uniformly in $r$. Let $f$ be the classical solution to \eqref{feco1} on the maximal interval $[0, \tau(\fin{f})),$ constructed in Theorem \ref{fcth1}, and let with  $f^{(r)} = f|_{X_{0,m,r}}$. Introducing the error $e_r(t,x) = f^{(r)}(t,x)-f_r(t,x)$, we have  $e_r(t,x)=0$ for $x>r$ and $e(0,x) =0$.  Using the fact that both $f$ (see \cite[Theorem 2.2]{BaLa12a}) and $f_r$ satisfy \eqref{PhPr001a} pointwise, for $0<x<r$ we have
\begin{align}
\p_t e_r(t,x) &= -a(x)e_r(t,x) + \cl{0}{\infty}a(x) b(x,y) e_r(t,y)\md x\label{e1}\\& -f^{(r)}(t,x)\cl{0}{\infty}k(x,y)f(t,y)\md y + \frac{1}{2}\cl{0}{x}k(x-y,y)f^{(r)}(t,x-y)f^{(r)}(t,y)\md y\nn\\
&+f_r(t,x)\cl{0}{r-x}k_r(x,y)f_r(t,y)\md y - \frac{1}{2}\cl{0}{x}k_r(x-y,y)f_r(t,x-y)f_r(t,y)\md y.\nn
\end{align} Then, we transform the coagulation part as follows
\begin{align*}
&\phantom{x} -f^{(r)}(t,x)\cl{0}{r-x}k(x,y)f^{(r)}(t,y)\md y -f^{(r)}(t,x)\cl{r-x}{\infty}k(x,y)f(t,y)\md y\nn\\& \phantom{x}+ \frac{1}{2}\cl{0}{x}k(x-y,y)f^{(r)}(t,x-y)f^{(r)}(t,y)\md y\nn\\
&\phantom{x}+f_r(t,x)\cl{0}{r-x}k(x,y)f_r(t,y)\md y - \frac{1}{2}\cl{0}{x}k(x-y,y)f_r(t,x-y)f_r(t,y)\md y\nn\\
   &= -\cl{0}{r-x}k(x,y)(f^{(r)}(t,x)f^{(r)}(t,y)-f_r(t,x)f_r(t,y))\md y\nn\\&\phantom{x} + \frac{1}{2}\cl{0}{x}k(x-y,y)(f^{(r)}(t,x-y)f^{(r)}(t,y)-f_r(t,x-y)f_r(t,y))\md y \nn\\ &\phantom{x}-f^{(r)}(t,x)\cl{r-x}{\infty}k(x,y)f(t,y)\md y,\end{align*}
where, for $0\leq x\leq r,$  $(x-y) + y = x\leq r$ so that $k_r(x-y,y) = k(x-y,y)$. Next we write
$$
f^{(r)}(t,x)f^{(r)}(t,y)-f_r(t,x)f_r(t,y) =  f^{(r)}(t,y)e_r(t,x) + e_r(t,y)f_r(t,x),
$$
and hence \eqref{e1} takes the form
\begin{align}
\p_t e_r(t,x) &= -a(x)e_r(t,x) + \cl{0}{\infty}a(x) b(x,y) e_r(t,y)\md x\label{e2}\\& \phantom{x}-\cl{0}{r-x}k(x,y)(f^{(r)}(t,y)e_r(t,x) + e_r(t,y)f_r(t,x))\md y\nn\\& \phantom{x}+ \frac{1}{2}\cl{0}{x}k(x-y,y)(f^{(r)}(t,y)e_r(t,x-y) + e_r(t,y)f_r(t,x-y))\md y \nn\\ &\phantom{x}-f^{(r)}(t,x)\cl{r-x}{\infty}k(x,y)f(t,y)\md y =: F_{0,m} e_r(t,x) + E_1 e_r(t,x) + E_2(f^{(r)},f)(t,x).\nn
\end{align}
Next, by  \eqref{cest2},
\begin{align*}
\|E_1 e_r(t)\|_{0,m}&=\cl{0}{r}w_m(x)\left|-\cl{0}{r-x}k(x,y)(f^{(r)}(t,y)e_r(t,x) + e_r(t,y)f_r(t,x))\md y\right. \\
&\phantom{x} + \left. \frac{1}{2}\cl{0}{x}k(x-y,y)(f^{(r)}(t,y)e_r(t,x-y) + e_r(t,y)f_r(t,x-y))\md y \right|\md x \\
&\leq L_0\|e_r(t)\|_{0,m}^{(\alpha)}(\|f_r(t)\|_{0,m}^{(\alpha)} + \|f^{(r)}(t)\|_{0,m}^{(\alpha)}) \leq L \|e_r(t)\|_{0,m}^{(\alpha)},
\end{align*}
where $L$ is a constant independent of $r$. Further,
\begin{align*}
\|E_2(f^{(r)},f)(t)\|_{0,m}&=\cl{0}{r} w_m(x) f^{(r)}(t,x)\left(\,\cl{r-x}{\infty} k(x,y)f(t,y)\md y \right)\md x\\& \leq \cl{0}{r}\frac{ w_m(x) f^{(r)}(t,x)}{1+(r-x)^m}\left(\,\cl{r-x}{\infty} k(x,y)f(t,y)w_m(y)\md y \right)\md x\\
&\leq K\|f(t)\|_{0,m}^{(\alpha)} \cl{0}{r}\frac{ w_m(x)(1+a(x))^\alpha f^{(r)}(t,x)}{1+(r-x)^m}\md x\\
&\leq K\|f(t)\|_{0,m}^{(\alpha)}\left(\frac{ \|f(t)\|_{0,m}^{(\alpha)}}{1+ \left(\frac{r}{2}\right)^m} + \|f(t)-f^{(\frac{r}{2})}(t)\|_{0,m}^{(\alpha)}\right).
\end{align*}
Then, using the integral formulation of \eqref{e2}, $e_r(0)=0$ and identifying $G_{F^{(r)}_{0,m}}(t) = G_{F_{0,m}}(t)$ on $X_{0,m,r}$, we have
\begin{equation}
e_r(t) = \cl{0}{t} G_{F_{0,m}}(t-s)E_1 e_r(s)\md s + \cl{0}{t} G_{F_{0,m}}(t-s)E_2(f^{(r)},f)(s)\md s.
\end{equation}
Using the analyticity of \sem{G_{F_{0,m}}}, \eqref{eq2.4a} and the estimates above, we have, for any fixed $t<\tau(\fin{f})$,
\begin{align*}
\|e_r(t)\|_{0,m}^{(\alpha)} \leq c_1 \cl{0}{t} \frac{\|e_r(s)\|_{0,m}^{(\alpha)}\md s}{(t-s)^\alpha} + c_2\left(\frac{1}{1+ \left(\frac{r}{2}\right)^m} + \sup\limits_{0\leq s\leq t}\|f(s)-f^{(\frac{r}{2})}(s)\|_{0,m}^{(\alpha)}\right),
\end{align*}
where $c_1,c_2$ are uniform in $r$ and $t \in [0,\tau(f^{in}))$. Using now the Gronwall-Henry inequality in the form of \cite[Lemma 3.2]{Banasiak2019}, we obtain
\begin{equation}
\|e_r(t)\|_{0,m}^{(\alpha)} \leq c_3\left(\frac{1}{1+ \left(\frac{r}{2}\right)^m} + \sup\limits_{0\leq s\leq t}\|f(s)-f^{(\frac{r}{2})}(s)\|_{0,m}^{(\alpha)}\right)
\label{finerr}
\end{equation}
for some constant $c_3$ independent of $r$. We observe that for any $r_n\to \infty$, $(\|f(t)-f^{(\frac{r_n}{2})}(t)\|_{0,m}^{(\alpha)})_{n\in \mbb N}$ is a monotone sequence of continuous functions converging to 0 (which is also a continuous function) and hence the convergence is uniform by Dini's theorem. Thus
$$
\lim\limits_{r \to \infty} \|e_r(t)\|_{0,m}^{(\alpha)} = 0
$$
uniformly in $t$ on any interval $[0,T]\subset [0,\tau(\fin f))$. Then also the approximate solutions $(f_r)_{r>0},$ extended by 0 to $\mbb R_+$, converge  to $f$ in $C([0,T],X_{0,m}^{(\alpha)})$. \hfill \qed

\section{Global solvability} In this section we prove the following theorem
\begin{theorem}
Assume that $a$ and $b$  satisfy \eqref{aloc}, \eqref{PhPr004}, \eqref{ass2} and \eqref{goodchar}, and  let $m > \max\{1,l\}$. If the coagulation kernel $k$ satisfies (\ref{kass2}), then, for each $f^{in}\in X_{0,m,+}^{(\alpha)}$, the corresponding local nonnegative classical solution \eqref{fprop}  is global in time.
\end{theorem}
  \proof Let us fix some $m_0$ for which the assumptions of Theorem \ref{fcth1} are satisfied. By a standard argument, we can assume that $f$ is defined on its maximal forward interval of existence $ [0, \tau(f^{in}))$. By \cite[Proposition 7.1.8]{Lun} (and the comment below it),  if $\tau(f^{in})<\infty$, then  $t\mapsto \|
f(t)\|_{[0,m_0]}^{(\alpha)}$ is unbounded as $t\to \tau(f^{in})$. Thus, to prove that $f$ is globally defined, we need to show that  $t\mapsto \| f(t)\|_{[0,m_0]}^{(\alpha)}$ is \textit{a priori} bounded on $[0,\tau(f^{in}))$.   We use the following two observations. First, for $0 \leq m_1 \leq m_2,$  $X_{0,m_2}$ is continuously and densely imbedded in $X_{0,m_1}$ and hence the boundedness of  $t\mapsto \Vert f(t) \Vert_{[0,m_2]}$  implies the same for $t\mapsto \Vert f(t)\Vert_{[0,m_1]}$ for all $ m_1 \in [0,m_2]$. Second, if Theorem \ref{fcth1} holds for a specific  $m_0$,  then it is also valid  in the scale of spaces $X_{0,m}$ with $m > \max\{1,l\},$ hence we can always choose an integer $m \geq \max\{2, m_0\}$ for which  Theorem~\ref{fcth1} holds.

Next, we need to establish several inequalities valid in any space $X_{0,i,+}.$  By  \eqref{moments2}, (\ref{momestb})  and assumption (\ref{kass2}), we can deduce that for each  $ i> 1$ and $f \in X_{0,i,+}$,
\begin{eqnarray}
&&\int_{0}^{\infty}x^i\mc Cf(x)\md x =
\frac{1}{2}\int_{0}^{\infty}\int_{0}^{\infty}((x+y)^i-x^i-y^i)k(x,y)f(x)f(y)\md x\md y\nn\\
&&\leq K_i(\|f\|^{(\alpha)}_{[i-1]}\|f\|_{[1]} + \|f\|_{[i-1]}\|f\|_{[1]}^{(\alpha)}),
\label{Cmom1}
\end{eqnarray}
where $K_i$ is a positive constant. For the case $i=1$ we have
\[
\int_{0}^{\infty}x\mc Cf(x)\md x = 0.
 \]
Turning now to the linear terms in \eqref{feco1},  we recall from  Theorem \ref{thnewchar0}, item 1.,  that, if $N_{m_0}(x)/x^{m_0} \geq \delta_{m_0}$ holds for some $m_0 > 1,$ $\delta_{m_0}$ and $x\geq 1$, then  there is $\delta_i>0$ such that $N_i(x)/x^i \geq \delta_i>0$ for any $i > 1$ and $x\geq 1$. Hence, for $f \in D(A_{0,i})_+$,
  \begin{align}
  &\cl{0}{\infty}(\mc Af(x) +\mc Bf(x))x^i\md x
=
-\cl{0}{\infty}N_i(x)a(x)f(x)\md x \nn\\&=  -\cl{0}{1}a(x) N_i(x)f(x)\md x - \cl{1}{\infty} (a(x)+\omega)f(x) x^i \frac{N_i(x)}{x^i}\md x + \omega \cl{1}{\infty} f(x) N_i(x) \md x\nn\\
&\leq  -\cl{0}{1}a(x) N_i(x)f(x)\md x - \delta_i\cl{1}{\infty} (a(x)+\omega)f(x) x^i \md x + \omega \cl{1}{\infty} f(x) N_i(x) \md x\nn\\
&= - \delta_i \|f\|_{[i]}^{(1)} -\cl{0}{1}a(x) N_i(x)f(x)\md x + \delta_i \cl{0}{1}(a(x)+\omega) x^if(x)\md x + \omega\cl{1}{\infty} f(x) N_i(x) \md x\nn\\
&\leq - \delta_i \|f\|_{[i]}^{(1)} + \omega_1 \|f\|_{[i]},\label{fragest}
\end{align}
 where $\omega_1 = \delta_i\text{ess}\sup_{0\leq x\leq 1} a(x) + \omega(1+\delta_i)$. As for the coagulation term, for $i=1$ we have
\[
\int_{0}^{\infty}x(\mc Af(x) +\mc Bf(x))\md x = 0.
\]
 If we take $\fin{f}$ with bounded support in $[0,\infty)$, then $\fin{f} \in X^{(\alpha)}_{0,i,+}$ and, if additionally $i>\max\{1,l\}$, then the corresponding solution $(0,\tau(f^{in}))\ni t\mapsto f(t)$ is differentiable in any such $X^{\alpha}_{0,i}$ and thus in any $X_{i}, i\geq 0,$ or, in other words, any moment of the solution is differentiable. First, let us consider an integer $i \geq 2$. Then, from \eqref{Cmom1} and \eqref{fragest},  we have
\begin{eqnarray}
 \frac{d}{dt} \|f(t)\|_{[i]} &\leq& \omega_1 \|f(t)\|_{[i]} - \delta_i \|f(t)\|_{[i]}^{(1)} \nn \\
 && \quad + K_i(\|f(t)\|^{(\alpha)}_{[i-1]}\|f(t)\|_{[1]} + \|f(t)\|_{[i-1]}\|f(t)\|_{[1]}^{(\alpha)}).\label{firstmom}
\end{eqnarray}
To simplify \eqref{firstmom}, we use  the following auxiliary inequalities.  For  $i \geq 2$ and  $1\leq r \leq i-1,$ we apply the H\"{o}lder's inequality with $p=1/\alpha$ and $q =1/(1-\alpha)$ to  obtain
\begin{align}
\|f\|_{[r]}^{(\alpha)} &= \int_0^\infty x^r a_\omega^\alpha (x) f(x) \md x = \int_0^1 x^r a_\omega^\alpha (x) f(x) \md x + \int_1^\infty x^r a_\omega^\alpha (x) f(x) \md x\nn\\
&\leq c_a\int_0^1 x f(x) \md x  + \int_1^\infty x^{(i-1)/q}f^{1/q}(x)x^{(qr-i+1)/q} a_\omega^\alpha (x) f^{1/p}(x) \md x\nn\\
&\leq c_a\| f\|_{[1]}  + \left(\int_0^\infty x^{i-1}f(x)\md x\right)^{1-\alpha}\left(\int_1^\infty x^{(r-(i-1)(1-\alpha))/\alpha} a_\omega(x) f(x) \md x\right)^{\alpha}\nn\\
&\leq c_a\| f\|_{[1]}  + \|f\|_{[i-1]}^{1-\alpha}\left(\|f\|^{(1)}_{[i]}\right)^{\alpha}. \label{wl866}
\end{align}
Note that the above derivation of \eqref{wl866} uses the fact that $(r-(i-1)(1-\alpha))/\alpha \leq i-1 < i$ for $\alpha\in (0,1)$ and $r\leq i-1,$ and hence
$$
x^{(r-(i-1)(1-\alpha))/\alpha} \leq x^i, \qquad x \in [1,\infty).
$$
Young's inequality, with $p = 1/\alpha$ and $q = 1/(1-\alpha)$, then leads to
\begin{align}
\|f\|^{(\alpha)}_{[i-1]}\|f\|_{[1]} &\leq c_a\|f\|^2_{[1]} + \|f\|_{[1]}\|f\|_{[i-1]}^{1-\alpha}\left(\|f\|^{(1)}_{[i]}\right)^{\alpha}\nn\\
&\leq c_a\|f\|^2_{[1]} + \|f\|_{[1]}\left((1-\alpha)\e^{1/(\alpha -1)}\|f\|_{[i-1]} + \alpha\e^{1/\alpha}\|f\|^{(1)}_{[i]}\right)
\end{align}
and
\begin{align}
\|f\|_{[i-1]}\|f\|_{[1]}^{(\alpha)} &\leq c_a\|f\|_{[1]} \|f\|_{[i-1]}+ \|f\|_{[i-1]}^{2-\alpha}\left(\|f\|^{(1)}_{[i]}\right)^{\alpha}\nn\\
&\leq c_a\|f\|_{[1]} \|f\|_{[i-1]} + \left((1-\alpha)\e^{1/(\alpha-1)}\|f\|_{[i-1]}^{(2-\alpha)/(1-\alpha)} + \alpha\e^{1/\alpha}\|f\|^{(1)}_{[i]}\right).
\end{align}
We now apply these inequalities to the solution $t\mapsto f(t)$.  Since $\|f(t)\|_{[1]} = \|f^{in}\|_{[1]}$ is constant on $[0, \tau(f^{in}))$, by choosing $\e$ so that $\alpha\e^{1/\alpha} K_i(\|f\|_{[1]} +1)\leq \delta_i$, we see that there are positive constants $D_{0,i}, D_{1,i}, D_{2,i}, D_{3,i}$ such that (\ref{firstmom}) can be written as
\begin{equation}
 \frac{d}{dt} \|f(t)\|_{[i]} \leq D_{0,i}+ D_{1,i} \|f(t)\|_{[i]} + D_{2,i} \|f(t)\|_{[i-1]} + D_{3,i} \|f(t)\|^{(2-\alpha)/(1-\alpha)}_{[i-1]}.   \label{firstmom1}
\end{equation}
In particular, for $i = 2$ we obtain
\begin{equation}
 \frac{d}{dt} \|f(t)\|_{[2]} \leq D_{0,2}+ D_{1,2} \|f(t)\|_{[2]} + D_{2,2} \|f^{in}\|_{[1]} + D_{3,2} \|f^{in}\|^{(2-\alpha)/(1-\alpha)}_{[1]},   \label{firstmom2}
\end{equation}
and thus $t \mapsto \|f(t)\|_{[2]}$ is bounded on $[0, \tau(f^{in}))$. Then we can use \eqref{firstmom1} to proceed inductively to establish the  boundedness of $t \mapsto \|f(t)\|_{[i]}$ for all integer $i$. Further, since for any $i>1$ we have $x^i \leq x$ for $x \in [0,1]$ and $x^i \leq x^{\lfloor i\rfloor +1}$
 $$
 \|f\|_{[i]} \leq \|f\|_{[1]} + \|f\|_{[\lfloor i\rfloor +1]},
 $$
 we find that all moments of the solution of order $i\geq 1$ are bounded on the maximal interval of its existence.

 In the next step we show that also the moment of order 0 is bounded. In the proof we use the ideas of \cite[Theorem 8.2.23]{BLL}  but while there the estimates were  carried out for compactly supported approximating solutions, as in Theorem \ref{thappr}, here we will work with classical solutions defined for $x\in \mbb R_+$ but only  for $t\in [0, \tau(f^{in}))$.

Let us fix an integer $i>\max\{1,l\}.$ For the fragmentation term we have, as in  \eqref{fragest},
 \begin{equation}
 \cl{0}{\infty} [\mc F f](t,x)x^i\md x \leq -\delta_i\cl{1}{\infty} a(x)f(t,x) x^i\md x.
 \label{Fest}
 \end{equation}
 Let us define
 $$
 \Phi(t) : = \|f(t)\|_{[i]} + \delta_i\cl{0}{t}\cl{1}{\infty} a(x) f(s,x)x^i\md x\md s.
 $$
 We observe that, by induction,  \eqref{firstmom1}  can be written as \begin{equation}
 \frac{d \|f(t)\|_{[i]}}{dt}  \leq D_{0,i}+ D_{1,i} \|f(t)\|_{[i]} +  \Theta(t),  \label{firstmom1a}
\end{equation}
 where $\Theta(t)$ is bounded on $[0,\tau(\fin f))$. Then
 \begin{align*}
 \frac{d\Phi(t)}{dt} &= \frac{d \|f(t)\|_{[i]}}{dt}+ \delta_i\cl{1}{\infty} a(x) f(t,x)x^i\md x  \leq D_{0,i} + D_{1,i}\Phi(t) + \Theta(t)
 \end{align*}
and, integrating,
\begin{align*}
\Phi(t) &\leq e^{D_{1,i} t}\left(\Phi(0) + \frac{D_{0,i}}{D_{1,i}}(1-e^{-D_{1,i}t}) + \cl{0}{t}\Theta(s)e^{-D_{1,i} s}\md s\right)
\end{align*}
and we see that neither $\Phi,$ nor
\begin{equation}
t\mapsto \cl{0}{t}\cl{1}{\infty} a(x) f(s,x)x^i\md x\md s
\label{Pt}
\end{equation}
can blow up at $t=\tau(f^{in})$.
Let us define
$$
P(t):= \cl{1}{\infty} a(x) f(s,x)w_i(x)\md x\md s.
$$
Using the fact that
$$
\cl{0}{\infty} \mc C f(t,x)\md x \leq 0$$
and, by \eqref{4},
\begin{align}
\cl{0}{\infty} \mc F f(t,x)\md x &\leq \cl{0}{\infty} (n_0(y)-1)a(y) f(t,y)\md y \leq 2b_0\cl{0}{\infty} a(y) f(t,y) w_i(y)\md y \nn\\
&\leq a_1\cl{0}{1}  f(t,y) \md y + 2b_0 P(t),
\label{eqP1}
\end{align}
on $[0,\tau(\fin f)),$ where $a_1 = 2b_0 \text{ess}\sup_{y\in [0,1]}a(y)w_i(y) $, for the zeroth moment we have
$$
\frac{d}{dt}\|f(t)\|_{[0]} \leq a_1\|f(t)\|_{[0]} + 2b_0 P(t)
$$
and hence
$$
\|f(t)\|_{[0]} \leq e^{a_1 t}\left(\|\fin f\|_{[0]} + 2b_0\cl{0}{t} P(s)\md s\right).
$$
Now, using \eqref{Pt},
\begin{align}
 P(t)&=\cl{0}{t}\cl{1}{\infty} a(x) f(s,x)w_i(x)\md x\md s \leq 2\cl{0}{t}\cl{1}{\infty} a(x) f(s,x)x^i\md x\md s.
\label{Pt1}
\end{align}
is bounded on $[0,\tau(\fin f))$ and hence $\|f(t)\|_{[0]}$ is bounded there as well.

To complete the proof, let $m_0>\max\{1,l\}$ be arbitrary. From the previous part of the proof, for $\fin{f}$ with bounded support, the norm of the corresponding solution, $t \mapsto \|f(t)\|_{[0,m_0]},$ remains bounded on $[0,\tau(\fin{f}))$. We again use the properties of the analytic semigroup \sem{F_{0,m_0}}. The local solution $f$ satisfies the integral equation
\begin{equation}
f(t) = G_{F_{0,m_0}}(t)f^{in} + \cl{0}{t}G_{F_{0,m_0}}(t-s)[\mc Cf](s)\md s
\end{equation}
and thus, using \eqref{eq2.4a} and \eqref{cest3},
\begin{align}
\|f(t)\|_{[0,m_0]}^{(\alpha)} &\leq C_1 \|f^{in}\|_{[0,m_0]}^{(\alpha)} + C_2\cl{0}{t}  \frac{\|\mc C f(s)\|_{[0,m_0]}}{(t-s)^\alpha} \md s \nn\\
&\leq C_1 \|f^{in}\|_{[0,m_0]}^{(\alpha)} + 2^{m_0+1}KC_2\cl{0}{t}  \frac{\|f(s)\|_{[0,m_0]}^{(\alpha)}\|f(s)\|_{[0,m_0]}}{(t-s)^\alpha} \md s\nn\\&\leq C_3 + C_4\cl{0}{t}  \frac{\|f(s)\|_{[0,m_0]}^{(\alpha)}}{(t-s)^\alpha} \md s,
\end{align}
where $C_3$ and $C_4$ are independent of time on $[0,\tau(f^{in}))$ on account of the boundedness of $t\mapsto \|f(t)\|_{[0,m_0]}$. Thus, using Gronwall--Henry inequality, for some constant $C_5$ independent of $t$,
\begin{equation}
\|f(t)\|_{[0,m_0]}^{(\alpha)} \leq C_5, \quad t\in [0,\tau(f^{in}))
\end{equation}
and hence $t \mapsto f(t)$ is a global classical solution in any $X_{0,m_0}^{\alpha}$ for which the assumptions of the theorem hold.

 To prove  the global existence of solutions emanating from any initial condition $f^{in}\in X_{0,m_0,+}^{\alpha}$ (and also to mild solutions) we observe that since the space of functions with bounded support is dense in $X_{0, m_0}^{\alpha}$ (respectively, $X_{0,m_0}$), a finite-time blow-up of such a solution would contradict the theorem on the continuous dependence of solutions on the initial data (which, in this case, follows from the Gronwall--Henry inequality, see \cite[Theorem 7.1.2]{Lun}) along the lines of the proof of \cite[Theorem 8.1.1]{BLL}. \qed


\begin{thebibliography}{99}
\bibitem{AizBak}
M.~Aizenman and T.~A. Bak.
\newblock Convergence to equilibrium in a system of reacting polymers.
\newblock {\em Comm. Math. Phys.}, 65(3):203--230, 1979.
\bibitem{AR}
W.~Arendt and A.~Rhandi.
\newblock Perturbation of positive semigroups.
\newblock {\em Arch. Math. (Basel)}, 56(2):107--119, 1991.
\bibitem{BaCa90}
J.~M. Ball and J.~Carr.
\newblock The discrete coagulation-fragmentation equations: existence,
  uniqueness, and density conservation.
\newblock {\em J. Statist. Phys.}, 61(1-2):203--234, 1990.
\bibitem{BaAr}
J.~Banasiak and L.~Arlotti.
\newblock {\em Perturbations of positive semigroups with applications}.
\newblock Springer Monographs in Mathematics. Springer-Verlag London, Ltd.,
  London, 2006.
  \bibitem{BLL}
J.~Banasiak, W.~Lamb and P.~Lauren\c{c}ot.
\newblock {\em Analytic Methods for Coagulation-Fragmentation Models, Volume I \& II}.
\newblock Chapman \& Hall/CRC Monographs and Research Notes in Mathematics,
  CRC Press, Boca Raton, 2019 (in print).
  \bibitem{Banasiak2019}
J.~{Banasiak}, L.~O. {Joel}, and S.~{Shindin}.
\newblock The discrete unbounded coagulation-fragmentation equation with
  growth, decay and sedimentation.
\newblock ArXiv e-prints, arXiv:1809.00046, 2018.
\bibitem{BaLa12a}
J.~Banasiak and W.~Lamb.
\newblock Analytic fragmentation semigroups and continuous
  coagulation-fragmentation equations with unbounded rates.
\newblock {\em J. Math. Anal. Appl.}, 391(1):312--322, 2012.
\bibitem{BLL13}
J.~Banasiak, W.~Lamb and M.~Langer.
\newblock Strong fragmentation and coagulation with power-law rates.
\newblock {\em J. Engrg. Math.}, 82:199--215, 2013.
\bibitem{becker}
R.~Becker and W.~D{\"o}ring.
\newblock Kinetische {B}ehandlung der {K}eimbildung in {\"u}bers{\"a}ttigten
  {D}{\"a}mpfen.
\newblock {\em Annalen der Physik}, 416(8):719--752, 1935.
\bibitem{bergh1976}
J.~Bergh and J.~L{\"o}fstr{\"o}m.
\newblock \emph{Interpolation spaces: an introduction}.
\newblock Springer-Verlag, Berlin-New York, 1976.
\bibitem{Bert06}
J.~Bertoin.
\newblock {\em Random fragmentation and coagulation processes}, volume 102 of
  {\em Cambridge Studies in Advanced Mathematics}.
\newblock Cambridge University Press, Cambridge, 2006.
\bibitem{blatz1945}
P.~J. Blatz and A.~V. Tobolsky.
\newblock Note on the kinetics of systems manifesting simultaneous
  polymerization-depolymerization phenomena.
\newblock {\em The journal of Physical Chemistry}, 49(2):77--80, 1945.
\bibitem{DuSt96b}
P.~B. Dubovski{\u\i} and I.~W. Stewart.
\newblock Existence, uniqueness and mass conservation for the
  coagulation-fragmentation equation.
\newblock {\em Math. Methods Appl. Sci.}, 19(7):571--591, 1996.
  \bibitem{EN}
K.-J. Engel and R.~Nagel.
\newblock {\em One-parameter semigroups for linear evolution equations}, volume
  194 of {\em Graduate Texts in Mathematics}.
\newblock Springer-Verlag, New York, 2000.
\bibitem{EMP02}
M.~Escobedo, S.~Mischler, and B.~Perthame.
\newblock Gelation in coagulation and fragmentation models.
\newblock {\em Comm. Math. Phys.}, 231(1):157--188, 2002.
\bibitem{ELMP03}
M.~Escobedo, P.~Lauren{\c{c}}ot, S.~Mischler, and B.~Perthame.
\newblock Gelation and mass conservation in coagulation-fragmentation models.
\newblock {\em J. Differential Equations}, 195(1):143--174, 2003.
\bibitem{LaMi02b}
{\relax Ph}.~Lauren{\c{c}}ot and S.~Mischler.
\newblock From the discrete to the continuous coagulation-fragmentation
  equations.
\newblock {\em Proc. Roy. Soc. Edinburgh Sect. A}, 132(5):1219--1248, 2002.
\bibitem{Lun}
A.~Lunardi.
\newblock {\em Analytic semigroups and optimal regularity in parabolic
  problems}, volume~16 of {\em Progress in Nonlinear Differential Equations and
  their Applications}.
\newblock Birkh\"auser Verlag, Basel, 1995.
\bibitem{McZi87}
E.~D. McGrady and R.~M. Ziff.
\newblock ``{S}hattering'' transition in fragmentation.
\newblock {\em Phys. Rev. Lett.}, 58(9):892--895, 1987.
\bibitem{Melz57b}
Z.~A. Melzak.
\newblock A scalar transport equation.
\newblock {\em Trans. Amer. Math. Soc.}, 85:547--560, 1957.
\bibitem{muller1928allgemeinen}
H.~M{\"u}ller.
\newblock Zur allgemeinen {T}heorie der raschen {K}oagulation.
\newblock {\em Fortschrittsberichte {\"u}ber Kolloide und Polymere},
  27(6):223--250, 1928.
\bibitem{Smo16}
M.~v. Smoluchowski.
\newblock Drei {V}ortrage \"{u}ber {D}iffusion, {B}rownsche {B}ewegung und
  {K}oagulation von {K}olloidteilchen.
\newblock {\em Zeitschrift f{\"u}r Physik}, 17:557--585, 1916.
\bibitem{Smo17}
M.~v. {Smoluchowski}.
\newblock Versuch einer mathematischen {T}heorie der {K}oagulationskinetik
  Kolloider {L}{\"o}sungen.
\newblock {\em Zeitschrift f{\"u}r Physikalische Chemie}, 92:129 -- 168, 1917.
\bibitem{Stew89}
I.~W. Stewart.
\newblock A global existence theorem for the general coagulation-fragmentation
  equation with unbounded kernels.
\newblock {\em Math. Methods Appl. Sci.}, 11(5):627--648, 1989.
\bibitem{Stew91}
I.~W. Stewart.
\newblock Density conservation for a coagulation equation.
\newblock {\em Z. Angew. Math. Phys.}, 42(5):746--756, 1991.
  \bibitem{vigi}
R.~D. Vigil and R.~M. Ziff.
\newblock On the scaling theory of two-component aggregation.
\newblock {\em Chemical Engineering Science}, 53(9):1725--1729, 1998.
\bibitem{Voi77}
J.~Voigt.
\newblock On the perturbation theory for strongly continuous semigroups.
\newblock {\em Math. Ann.}, 229(2):163--171, 1977.
\end{thebibliography}
\end{document}